\documentclass[12pt]{article}
\usepackage{latexsym}

\newtheorem{Definition}{Definition}[section]
\newtheorem{Proposition}{Proposition}
\newtheorem{Lemma}{Lemma}[section]

\newcommand {\cC}{\mbox{${\mathcal C}$}}
\newcommand {\cD}{\mbox{${\mathcal D}$}}

\newcommand {\cS}{\mbox{${\mathcal S}$}}

\newcommand {\id}{\mbox{${\mathrm{id}}$}}

\usepackage{url}

\begin{document}

\title{Category Free Category Theory and Its Philosophical Implications}

\author{Michael Heller\thanks{Correspondence address: ul.
Powsta\'nc\'ow Warszawy 13/94, 33-110 Tarn\'ow, Poland. E-mail:
mheller@wsd.tarnow.pl} \\
Copernicus Center for Interdisciplinary Studies \\
ul. S{\l}awkowska 17, 31-016 Cracow, Poland}

\date{\today}

\maketitle

\begin{abstract}
There exists a dispute in philosophy, going back at least to Leibniz, whether is it possible to view the world as a network of relations and relations between relations with the role of objects, between which these relations hold, entirely eliminated. Category theory seems to be the correct mathematical theory for clarifying conceptual possibilities in this respect. In this theory, objects acquire their identity either by definition, when in defining category we postulate the existence of objects, or formally by the existence of identity morphisms. We show that it is perfectly possible to get rid of the identity of objects by definition, but the formal identity of objects remains as an essential element of the theory. This can be achieved by defining category exclusively in terms of morphisms and identity morphisms (objectless, or object free, category) and, analogously, by defining category theory entirely in terms of functors and identity functors (categoryless, or category free, category theory). With objects and categories eliminated, we focus on the ``philosophy of arrows'' and the roles various identities play in it (identities as such, identities up to isomorphism, identities up to natural isomorphism...). This perspective elucides a contrast between ``set ontology'' and ``categorical ontology''.
\end{abstract}

\section{Introduction}
There exists a tendency in philosophy, going back at least to Leibniz, to view the universe as a network of relations and relations between relations with the role of objects, between which these relations hold, minimalysed or even reduced to null. If in Leibniz such a view should be possibly restricted only to space and time\footnote{It is debatable whether it can be extended to his monadology.}, thinkers were not lacking who espoused this standpoint in its full extent. It is interesting to notice that one of Leibniz's arguments against Newton's idea of absolute space drew on his outspoken Principle of the Identity of Indiscernibles which excludes the possibility for two things being distinct that would not be different with respect to at least one discernible property. Imagine now the world moved by some distance in absolute space. Two such ``configurations'' (before and after the world has been moved) would be indiscernible even for God. Therefore, the concept of absolute space should be ruled out as contradicting the Principle of the Identity of Indiscernibles. Essentially the same reasoning leads, according to Leibniz, to rejection of absolute time. The concepts of identity and that of relatedness of everything with everything are strictly interwoven with each other. To see the interconnectedness of these concepts, let us introduce the following terminology. There is a sense in which ``we can talk about \textit{this} object, whatever properties it may have, as opposed to \textit{that} object, whatever properties this second object may have''\cite[pp. 16-17]{Teller}. In this sense, we speak about the \textit{primitive thisness} of this object. If this is not the case, i.e. if the identity of the object depends on its properties related to other objects, we can talk about the \textit{induced} or \textit{contextual thisness} of this object. If the identity of every object is reduced only to the induced thisness, we can speak about the \textit{complete relatedness} of the system of objects. On the contrary, if the identity of every object is reduced only to its primitive thisness, we could speak about an ``absolute'' system in which relations between objects play no role. Of course, ``mixed situations'' are also possible.

Category theory seems to be the correct mathematical theory for clarifying conceptual possibilities in this domain (see, for instance, \cite{Rodin05}). It seems natural to identify objects of the preceding paragraph with objects of category theory, and relations of the preceding paragraph with morphisms of category theory, or objects with categories themselves and relations with functors between categories (in this paper we do not consider $n$-categories). Is there a possibility to get rid of objects entirely and, in this way, to implement the idea of the complete relatedness of the system? To pose correctly this question in the context of category theory, we must first look for the identity of objects. In category theory objects acquire their identity on two levels: first, by definition, when in defining category we postulate the existence of objects; second, formally by the existence of identity morphisms; we can speak about the ``identity by definition'' and the ``formal identity'', respectively. In the following, we show that it is perfectly possible to get rid of the identity of objects by definition, but the formal identity of objects remains as an essential element of the theory.

\section{Objectless Category Theory}
In the present section, we briefly present, following \cite[pp. 44-46]{Semadeni} (see also \cite[pp. 41-43]{Adamek}), a formulation of category theory in terms of morphisms without explicitly mentioning objects. As we shall see, this formulation is equivalent to the standard one.

\begin{Definition}
Let us axiomatically define a theory which we shall call an \emph{objectless} or {object free category theory}. In this theory, the only primitive concepts (besides the usual logical concepts and the equality concept) are: 

(I) $\alpha $ is a morphism, 

(II) the composition $\alpha \beta $ is defined and is equal to $\gamma $, 

The following axioms are assumed: 
\begin{enumerate}
	\item \emph{Associativity of compositions:} Let $\alpha , \beta , \gamma $ be morphisms. If the compositions $\beta \alpha $ and $\gamma \beta $ exist, then 
		\begin{itemize}
		\item
		the compositions $\gamma (\beta \alpha )$ and $(\gamma \beta ) \alpha $ exist and are equal;
		\item
		if $\gamma (\beta \alpha )$ exists, then $\gamma \beta $ exists, and if $(\gamma \beta )\alpha $ exists then  $\beta \alpha $ exists.
		\end{itemize}
	\item \emph{Existence of identities:} For every morphism $\alpha $ there exist morphisms $\iota$ and $\iota'$, called identities,  such that 
	\begin{itemize}
	\item
	$\beta \iota = \beta $ whenever $\beta \iota $ is defined (and analogously for $\iota' $),
	\item
	$\iota \gamma = \gamma $ whenever $\iota \gamma $ is defined (and analogously for $\iota' $).
	\item
	$\alpha \iota $ and $\iota' \alpha  $ are defined.
	\end{itemize}
	\end{enumerate}
\end{Definition}

\begin{Lemma}
Identities $\iota $ and $\iota' $ of axiom (2) are  uniquely determined by the morphism $\alpha $.
\end{Lemma}

\noindent
\textit{Proof.} Let us prove the uniqueness for $\iota $ (for $\iota'$ the proof goes analogously). Let $\iota_1 $ and $\iota_2$ be identities, and $\alpha \iota_1$ and $\alpha \iota_2 $ exist. Then $\alpha \iota_1 = \alpha $ and $(\alpha \iota_1 )\iota_2 = \alpha \iota_2$. From axiom (1) it follows that $\iota_1 \iota_2$ is determined. 
But $\iota_1 \iota_2$ exists if an only if $\iota_1 = \iota_2$. Indeed, let us assume that $\iota_1 \iota_2$ exist then $\iota_1 = \iota_1 \iota_2 = \iota_2$. And vice versa, assume that $\iota_1 = \iota_2$. Then from axiom (2) it follows that there exists an identity $\iota $ 
such that $\iota \iota_1$ exists, and hence is equal to $\iota$ (because $\iota_1$ is an identity). This, in turn, means that $(\iota \iota_1) \iota_2$ exists, because $(\iota \iota_1) \iota_2 = \iota \iota_2 = \iota \iota_1 = \iota$. Therefore, $\iota_1 \iota_2$ exists by Axiom 1.
$\Box $
\par

Let us denote by $d(\alpha )$ and $c(\alpha )$ identities that are uniquely determined by a morphism $\alpha $, i.e. such that the compositions $\alpha d(\alpha ) $ and $c(\alpha )\alpha $ exist (letters $d$ and $c$ come from ``domain'' and ``codomain'', respectively).

\begin{Lemma}
The composition $\beta \alpha $ exists if and only if $c(\alpha ) = d(\beta )$, and consequently,
$$
d(\beta \alpha ) = d(\alpha ) \;\;\; \mathit{and} \;\;\; c(\beta \alpha ) = c(\beta ).
$$
\end{Lemma}

\noindent
\textit{Proof.} Let $c(\alpha ) = d(\beta ) = \iota $, then $\beta \iota $ and $\iota \alpha $ exist. From axiom (1) it follows that there exists the composition $(\beta \iota)\alpha = \beta \alpha $. Let us now assume that $\beta \alpha $ exists, and let us put $\iota = c(\alpha )$. Then $\iota \alpha $ exists which implies that $\beta \alpha = \beta (\iota \alpha ) = (\beta \iota )\alpha$. Since $\beta \iota $ exists then $d(\beta ) = \iota $. $\Box $

\begin{Definition}
If for any two identities $\iota_1 $ and $\iota_2 $ the class
\[
\left\langle \iota_1 , \iota_2 \right\rangle = \{\alpha : d(\alpha) = \iota_1,\; c(\alpha ) = \iota_2 \},
\]
is a set then objectless category theory is called \textit{small}.
\end{Definition}

\begin{Definition}\label{ObjectlessCat}
Let us choose a class \cC \ of morphisms of the objectless category theory (i.e. \cC \ is a model of the objectless category theory), and let $\cC^0$ denote the class of all identities of \cC . If $\iota_1, \iota_2, \iota_3 \in \cC^0$, we define the composition 
\[
m^{\cC^0}_{\iota_1, \iota_2, \iota_3 }: \left\langle \iota_1, \iota_2 \right\rangle \times \left\langle \iota_2, \iota_3 \right\rangle \rightarrow \left\langle \iota_1, \iota_3\right\rangle 
\]
by $m^{\cC^0}(\alpha , \beta ) = \beta \alpha $. Class \cC \ is called \emph{objectless category}.
\end{Definition}

\begin{Proposition}
The objectless category definition (\ref{ObjectlessCat}) is equivalent to the standard definition of category.
\end{Proposition}

\noindent
\textit{Proof.} To prove the theorem it is enough to reformulate the standard category definition in the following way.
A category \cC \ consists of

(I) a collection $\cC^0$ of objects,

(II) for each $A, B \in \cC^0$, 

a collection $\left\langle A, B\right\rangle_{\cC^0}$ of morphisms from $A$ to $B$,

(III) for each  $A, B, C \in \cC^0$, if $\alpha \in \left\langle A, B\right\rangle_{\cC^0}$ and $\beta \in \left\langle B, C\right\rangle_{\cC^0}$, 
the composition 
$$
m^{\cC^0}: \left\langle A, B\right\rangle_{\cC^0} \times \left\langle B, C\right\rangle_{\cC^0} \rightarrow \left\langle A, C\right\rangle_{\cC^0}
$$
is defined by $m^{\cC^0}_{A,B,C}(\alpha , \beta )$.
The following axioms are assumed
\begin{enumerate}
\item
\textit{Associativity:} If $\alpha \in \left\langle A, B\right\rangle_{\cC^0} , \beta \in \left\langle B, C\right\rangle_{\cC^0} , \gamma \in \left\langle C, D\right\rangle_{\cC^0} $ then 
$$
\gamma (\beta \alpha ) = (\gamma \beta )\alpha .
$$
\item
\textit{Identities:} For every $B \in \cC^0$ there exists a morphism $\iota_B \in \left\langle B, B\right\rangle_{\cC^0}$ such that
$$
 \forall_{A\in \cC^0} \forall_{{\alpha \in \left\langle A, B\right\rangle }_{\cC^0}} \; \iota_B \alpha = \alpha,
$$ $$
\forall_{C\in \cC^0} \forall_{{\beta \in \left\langle B, C\right\rangle }_{\cC^0}} \; \beta \iota_B = \beta .
$$
\end{enumerate}
To see the equivalence of the two definitions it is enough to suitably replace in the above definition objects by their corresponding identities. $\Box$

This theorem creates three possibilities to look at the category theory: (1) the standard way, in terms of objects and morphisms, (2) the objectless way, in terms of morphisms only, (3) the hybrid way in which we take into account the existence of objects but, if necessary or useful, we regard them as identity morphisms. The hybrid way could be useful in some interpretative issues.
 
\section{Category Theory without Categories}
It is almost trivial to see that the above move towards getting rid of objects can be repeated on the level of categories themselves, i.e. towards formulating category theory entirely in terms of functors. The idea would be the same as above -- to substitute identity functors for categories. What we need first, is the objectless functor definition, but it is straightforward
\begin{Definition}
Let \cC \ and \cD \ be objectless categories (definition (\ref{ObjectlessCat})). The \emph{objectless covariant (contravariant) functor} $\Phi : \cC \rightarrow \cD $ is an assignment of morphisms of \cC \ to the morphisms of \cD \ in such a way that the compositions and the identities are preserved, i.e if $\alpha , \beta $ are morphisms of \cC \ then
$$\Phi (\alpha \beta ) = \Phi (\alpha ) \Phi (\beta ),$$
$$(\Phi (\alpha \beta ) = \Phi (\beta ) \Phi (\alpha )),$$
and if $\iota , \iota' \in \cC^0$ then $\Phi (\iota )$ and $\Phi (\iota')$ are identities associated with $\Phi (\alpha )$.
\end{Definition}

Since functors can be composed and there exist identity functors, it is straightforward to repeat the strategy of the preceding section and define ``categoryless, or category free, category theory''. The obviousness of the definition allows us to do this in a sketchy way.

\begin{Definition}
The only primitive concepts of the \emph{category theory without categories}  are:

(I) $\Phi $ is a functor,

(II) the composition of functors $\Phi \Psi $ is defined and is equal to $\Lambda $.

The axioms of the associativity of functor composition and of the existence of functor identities (analogous to axioms (1) and (2) of definition (1.1)) are assumed.
\end{Definition}

The results analogous to that of Lemma 2.1 and 2.2 are clearly true. Let $D(\Phi )$ and $C(\Phi )$ be identity functors such that the compositions $\Phi D(\Phi )$ and $C(\Phi ) \Phi $ exist. Then the composition $\Psi \Phi $ exists if and only if $C(\Phi ) = D(\Psi )$.

The above results are almost trivial from the formal point of view, but philosophically they are quite remarkable: the very concept of category is not indispensable in developing the theory itself. This remains in agreement with the practice of Eilenberg and Mac Lane who, in their seminal paper \cite{EilenbergMacLane45}, regarded categories as auxiliary constructs necessary only to ensure domains and codomains for morphisms.

\section{Contextual Thisness}

What has been achieved by introducing objectless categories and category theory without categories?\footnote{It is interesting to notice that William Lawvere in his seminal doctor thesis consequently eliminated objects with the help of identity morphisms \cite{Lawvere}.} Let us focus on objects (analogous things can be said about categories). As we have seen, in the objectless category theory, objects do not acquire their individuality by definition, but only formally through the identity morphisms. Instead postulating the existence of objects, one postulates the existence of morphisms which, when composed with other morphisms, change nothing. One cannot speak here about the primitive thisness in the sense described in Section 1, but only about the contextual thisness, and the context is now given by those morphisms with which a given identity morphism composes (to change nothing). I propose to call it \textit{compositional thisness}. This is the closest (but still faraway) from the primitive thisness (as it can be defined in set theory) that can be obtained in category theory. Let us notice that the principle of indiscernibility of identities, in its original form, does not hold for this kind of thisness -- two identities could be discernible through their compositions with different morphisms. We can, therefore, speak of \textit{discernibility through compositions}.

The elimination of objects is possible not only on the formal level, but also on the level of logical language. Let us consider an example. For some purposes it is important to define the concept of category in a first-order language. The most obvious of such languages is the one in which there are two sorts of variables: one ranging over objects and another over arrows. However, the language becomes simpler and more elegant if we abolish this dichotomy by adopting objectless category approach and replace objects by their corresponding identity arrows. This requires some toil but is certainly more in the spirit of category theory (see \cite[232-234]{Goldblatt}).

\section{Individuality of Categories}
In mathematics we are usually interested in the individuality of entities ``up to isomorphism''. In this sense, mathematical entities change their individuality depending on the context. For instance, in topology two homeomorphic spaces are regarded as the same (homeomorphism is a topological isomorphism), but in differential geometry two diffeomorphic differential manifolds are regarded as the same (diffeomorphism is an isomorphism of differential manifolds). The concept of isomorphism has its categorical counterpart. 

Let $X$ and $Y$ be objects in a category \cC . A morphism $f: X \rightarrow Y$ is said to be isomorphism if there is a morphism $g: Y \rightarrow X$ such that   $g \circ f = \iota_X$ and $f\circ g = \iota_Y$. Since $g$ is necessarily unique, we may write $g = f^{-1}$. As we can see, this is fully compatible with the identification of objects in terms of identity morphisms: the composition of the defining morphisms with their respective inverses must give the identity arrow. We are entitled to say that the context with respect to which the thisness of an object is determined is now broadened to the extent the nature of mathematical reasoning demands: thisness is not determined by the identity arrow of a given object alone, but with the help of identities corresponding to suitable compositions ($g \circ f$ and $f \circ g$). I shall call it \textit{thisness up to isomorphism}.

The isomorphism concept refers also to categories. A functor $F: \cC \rightarrow \cD $ from a category \cC \ to a category \cD \ is said to be an isomorphism, if it has the inverse $G: \cD \rightarrow \cC $ such that $G \circ F = \id_{\cC }$ and $F \circ G = \id_{\cD }$ (here we assume that $\id_{\cC } \in \cC^0$ and $\id_{\cD } \in \cD^0$	). In such a case, we say that the categories \cC \ and \cD \ are isomorphic, and we write $\cC \cong \cD $. It is not difficult to see that what identity morphisms are for objects, the identity functors are for categories. It turns out, however, that the concept of contextual thisness as applied to categories (their isomorphism) is too rigid from the point of view of category theory. Two categories can be equal (remaining in a bijective correspondence), they can be unequal but isomorphic, or they can be not even isomorphic but nevertheless equivalent from the categorical point of view. To define the latter notion, we must first introduce the concept of naturality.

Let us consider two functors $F$ and $G$ from a category \cC \ to a category \cD , $F: \cC \rightarrow \cD $ and  $G: \cC \rightarrow \cD $. A natural transformation from the functor $F$ to the functor $G$ is an assignment $\tau $ that associates every object $X$ of \cC \ with the arrow $\tau_X: F(X) \rightarrow G(X)$ of \cD \ in such a way that for any arrow $f: X \rightarrow Y$ of \cC \ one has $\tau_Y \circ F(f) = G(f) \circ \tau_X$ (to see this clearly draw the corresponding commutative diagram). The arrows $\tau_X$ are said to be the components of $\tau $. We could imagine the functors $F$ and $G$ as giving two pictures of the category \cC \ within the category \cD . If these two pictures are ``faithful'' to the original and between themselves the transformation between the functors $F$ and $G$ is natural. If every component $\tau_X $ of $\tau $ is an isomorphism in \cD , then $\tau $ is said to be a natural isomorphism. If this is the case, two above mentioned pictures of the category \cC \ within the category \cD \ are not only ``faithful'' but also ``exact''.

Let us now go back to the equivalence of categories from the categorical point of view. A functor $F: \cC \rightarrow \cD $ is said to be an equivalence of categories, if there is a functor $G: \cD \rightarrow \cC $ such that there exist natural isomorphisms $\tau : \id_{\cC } \cong G \circ F$ and $\sigma : \id_{\cD } \cong F \circ G$; in such a case, we say that the categories \cC \ and \cD \ are equivalent and write $\cC \simeq \cD $. In tis context, I shall speak about \textit{thisness up to equivalence}.

The notion of the equivalence of categories is clearly more tolerant than that of the isomorphism of categories, in the sense that isomorphic categories are equivalent but not vice versa. The relationship between these two notions can be made more transparent in the following way. A category \cC \ is said to be skeletal if its isomorphic objects are identical, i.e. $A \cong B$ iff $A = B$ where $A$ and $B$ are objects of \cC . In other words, for skeletal categories, ``isomorphic'' means the same as ``is equal'', that is to say primitive thisness means the same as contextual thisness. A full subcategory\footnote{If $\cC_0$ is a subcategory of \cC \ then every object of $\cC_0$  is an object of \cC , and a subcategory $\cC_0$ of a category \cC \ is a full subcategory of \cC, if \cC \ has no arrows other than the ones already present in $\cC_0$.} $\cC_0$ of a category \cC \ is called a skeleton of \cC \ if $\cC_0$ is skeletal and each object of \cC \ is isomorphic to the exactly one object of $\cC_0$. Every two skeletons of a category are isomorphic.

A skeleton $\cC_0$ of \cC \ is clearly equivalent to \cC , $\cC_0 \simeq \cC $, and if categories have the same skeleton, they are equivalent\footnote{All mathematical facts quoted in this section are standard in category theory; see for instance \cite[pp. 200-202]{Goldblatt}.}. 

Finally, let us quote a few examples. The category of all finite sets and all functions between them has the subcategory of all finite ordinal numbers as its skeleton. The category of all well ordered sets with suitable order preserving morphisms has the subcategory of all ordinal numbers as its skeleton. The category of all vector spaces over a fixed field $K$ and $K$-linear transformations has the subcategory of all $K^n$, for $n$ any cardinal number, as its skeleton.

It is interesting to notice that if the axiom of choice holds for a category \cC \ then it has a skeleton. This statement can be proved in the following way. The relation of isomorphism on the collection of objects of a category \cC \ is an equivalence relation. Let us choose one object from each equivalence class. The full subcategory of \cC \ formed from the chosen elements is a skeleton of \cC .

\section{The Arrow Philosophy}
Categories are complex entities. They are composed of objects and arrows subject to a few axiomatic exigencies which, although rather simple, lead to a rich structural variety. Since, as we have seen, objects can be eliminated, everything finally reduces to arrows. Often the ``philosophy of arrows'' is contrasted with the ``philosophy of elements'' (each of them with its primitive thisness) and it is said that they determine two different ``ontologies'' for mathematics. The set ontology is of reductionist type: in it everything is reduced to ``being an element of'', typified by the functor ``$\in $''. An element either belongs to a set, or does not belong to it which is a symptom of the two-valued classical logic. The categorical ontology is of referential type with arrows, as fundamental entities, indicating an action or reference. Being an element of an object is replaced by an arrow pointing to a given object. Such an arrow could be understood as denoting an action that consists of picking up something (an ``element'') in this object. Moreover, this picking up can be done with various degrees of ``intensity'' depending on the domain of a given arrow: with maximal intensity (for sure) if the domain is the terminal object of the category, and with lesser intensity (with lesser certainty) otherwise.

Something similar happens on the level of categories with functors playing the role of arrows. The above ``arrow strategy'' works also at this level. As we have seen, categories could be entirely eliminated with identity functors effectively replacing them.

At the bottom of the ``new perspective'', exploited in the present work, lies the concept of individuality. In the ``set theoretical ontology'', the individuality of a set element is given by its very ``thisness'', independently of its relationship to the environment. In the ``ontology'' determined by the strategy of arrows, one can also speak on thisnesses of some structures, but it is always contextual from the very beginning. As we have seen, it strongly depends on its environment. What we have called compositional thisness  depends on arrows with which the corresponding arrow composes to change nothing. Two antiparallel arrows (i.e. arrows pointing in two opposite directions) can in various ways interact with each other to produce identities. If they compose to produce just identity, they are said to define an isomorphism, and we have, in such circumstances, agreed  to speak about thisness up to isomorphism. If they compose to produce identity in a natural way, we have agreed to speak about thisness up to equivalence.

To this list of producing identities ``of various degrees'' we can add some more. Let us consider the pair of functors between categories \cC  \ and \cD :  
$F:\cD \rightarrow \cC $  and   $G:\cC \rightarrow \cD $, and two natural transformations
\[
\epsilon : FG \rightarrow \mathrm{Id}_{\cC },
\]
and
\[
\eta: \mathrm{Id}_{\cD } \rightarrow GF.
\]
With the help of these data one defines the \textit{adjunction} between categories \cC \ and \cD ; $F$ is called the left adjoint to $G$, and $G$ is called the right adjoint to $F$.\footnote{For the definition see, for example,  \cite[chapter 9]{Awodey} or in a very accessible way \cite[chapter 5]{Simmons}.} $\epsilon $ and $\eta $ are then called the counit and the unit of this adjunction, respectively. Here we have also an interaction between natural transformations producing suitable unit and counit. Adjunction can be regarded as a generalisation of the equivalence of categories \cite[p. 209]{Awodey}. This time unit and counit arise from interaction between functors rather than between categories. Adjunction turns out to be  ubiquitous in mathematics; it unveils intimate relationships between sometimes very faraway departments of mathematics. In this sense, we could speak on non-local context dependence.

This ``strategy of units'' tells us something  about the nature of the ``field of categories''; it is drastically unlike a collection of sets. For the sake of concreteness let us specialise to toposes, and let us consider a pair of adjoint functors $f^*: \cD \rightarrow \cC$ and $f_*: \cC \rightarrow \cD $ such that $f^*$ is left adjoint to $f_*$, and $f^*$ is also left exact\footnote{A left exact functor is a functor that preserves finite limits.}. Such functors define an \textit{admissible transformation} between two toposes, which in this context are called frameworks (the name being clearly motivated by an analogy with reference frame used in physics)\footnote{More precisely, the name ``framework'' is used for any topos having the natural numbers object.}. The shifting through an admissible transformation from one framework to another framework can lead to a shift in understanding a mathematical concept. Here is an example. Let us consider  the concept of real-valued continuous function on a topological space $X$. Such a function can be interpreted, in the topos \cS \ of constant sets, as a continuously varying real number. Let us now shift, via the admissible transformation, from \cS \ to the topos $Shv(X)$ of sheaves over $X$. The concept of real function in \cS \ shifts to the concept of real number as it is interpreted in $Shv(X)$. J. L. Bell \cite{Bell} describes the situation in the following way:
\begin{quote}
...shifting to $Shv(X)$ from \cS \ essentially amounts to placing oneself in a framework which is, so to speak, itself `moving along' with the variation over $X$ of the given variable real numbers. This causes the variation of any variable real number not to be `noticed' in $Shv(X)$; it is accordingly there regarded as being constant real number.
\end{quote}

This is only the illustration of the main Bell's idea that ``the absolute universe of sets [should] be relinquished in favour of a plurality of local mathematical frameworks''. The interpretation of any mathematical concept is not fixed, but changes with the change of a local framework. It should be emphasised that units and counits (defining adjointness of admissible transformations) play here the key role. The concept of identity seems to lie at the core of the categorical architecture.

\end{document}